
\magnification\magstep1
\baselineskip15pt
\hfuzz12.2pt
\advance\hsize2.9pt

\newread\AUX\immediate\openin\AUX=\jobname.aux
\newcount\relFnno
\def\ref#1{\expandafter\edef\csname#1\endcsname}
\ifeof\AUX\immediate\write16{\jobname.aux gibt es nicht!}\else
\input \jobname.aux
\fi\immediate\closein\AUX

\newread\AUX\immediate\openin\AUX=\jobname.aux
\newcount\relFnno
\def\ref#1{\expandafter\edef\csname#1\endcsname}
\ifeof\AUX\immediate\write16{\jobname.aux gibt es nicht!}\else
\input \jobname.aux
\fi\immediate\closein\AUX



\def\ignore{\bgroup
\catcode`\;=0\catcode`\^^I=14\catcode`\^^J=14\catcode`\^^M=14
\catcode`\ =14\catcode`\!=14\catcode`\"=14\catcode`\#=14\catcode`\$=14
\catcode`\&=14\catcode`\'=14\catcode`\(=14\catcode`\)=14\catcode`\*=14
\catcode`+=14\catcode`\,=14\catcode`\-=14\catcode`\.=14\catcode`\/=14
\catcode`\0=14\catcode`\1=14\catcode`\2=14\catcode`\3=14\catcode`\4=14
\catcode`\5=14\catcode`\6=14\catcode`\7=14\catcode`\8=14\catcode`\9=14
\catcode`\:=14\catcode`\<=14\catcode`\==14\catcode`\>=14\catcode`\?=14
\catcode`\@=14\catcode`\A=14\catcode`\B=14\catcode`\C=14\catcode`\D=14
\catcode`\E=14\catcode`\F=14\catcode`\G=14\catcode`\H=14\catcode`\I=14
\catcode`\J=14\catcode`\K=14\catcode`\L=14\catcode`\M=14\catcode`\N=14
\catcode`\O=14\catcode`\P=14\catcode`\Q=14\catcode`\R=14\catcode`\S=14
\catcode`\T=14\catcode`\U=14\catcode`\V=14\catcode`\W=14\catcode`\X=14
\catcode`\Y=14\catcode`\Z=14\catcode`\[=14\catcode`\\=14\catcode`\]=14
\catcode`\^=14\catcode`\_=14\catcode`\`=14\catcode`\a=14\catcode`\b=14
\catcode`\c=14\catcode`\d=14\catcode`\e=14\catcode`\f=14\catcode`\g=14
\catcode`\h=14\catcode`\i=14\catcode`\j=14\catcode`\k=14\catcode`\l=14
\catcode`\m=14\catcode`\n=14\catcode`\o=14\catcode`\p=14\catcode`\q=14
\catcode`\r=14\catcode`\s=14\catcode`\t=14\catcode`\u=14\catcode`\v=14
\catcode`\w=14\catcode`\x=14\catcode`\y=14\catcode`\z=14\catcode`\{=14
\catcode`\|=14\catcode`\}=14\catcode`\~=14\catcode`\^^?=14
\Ignoriere}
\def\Ignoriere#1\;{\egroup}

\newcount\itemcount
\def\resetitem{\global\itemcount0}\resetitem
\newcount\Itemcount
\Itemcount0
\newcount\maxItemcount
\maxItemcount=0

\def\FILTER\fam\itfam\tenit#1){#1}

\def\Item#1{\global\advance\itemcount1
\edef\TEXT{{\it\romannumeral\itemcount)}}%
\ifx?#1?\relax\else
\ifnum#1>\maxItemcount\global\maxItemcount=#1\fi
\expandafter\ifx\csname I#1\endcsname\relax\else
\edef\testA{\csname I#1\endcsname}
\expandafter\expandafter\def\expandafter\testB\testA
\edef\testC{\expandafter\FILTER\testB}
\edef\testD{\csname0\testC0\endcsname}\fi
\edef\testE{\csname0\romannumeral\itemcount0\endcsname}
\ifx\testD\testE\relax\else
\immediate\write16{I#1 hat sich geaendert!}\fi
\expandwrite\AUX{\neverexpand\ref{I#1}{\TEXT}}\fi
\item{\ifx?#1?\relax\else\marginnote{I#1}\fi\TEXT}}

\def\today{\number\day.~\ifcase\month\or
  Januar\or Februar\or M{\"a}rz\or April\or Mai\or Juni\or
  Juli\or August\or September\or Oktober\or November\or Dezember\fi
  \space\number\year}
\font\sevenex=cmex7
\scriptfont3=\sevenex
\font\fiveex=cmex10 scaled 500
\scriptscriptfont3=\fiveex
\def\A{{\bf A}}
\def\G{{\bf G}}
\def\P{{\bf P}}

\def\phi{\varphi}
\def\epsilon{\varepsilon}
\def\theta{\vartheta}
\def\uauf{\lower1.7pt\hbox to 3pt{%
\vbox{\offinterlineskip
\hbox{\vbox to 8.5pt{\leaders\vrule width0.2pt\vfill}%
\kern-.3pt\hbox{\lams\char"76}\kern-0.3pt%
$\raise1pt\hbox{\lams\char"76}$}}\hfil}}

\def\title#1{\par
{\baselineskip1.5\baselineskip\rightskip0pt plus 5truecm
\leavevmode\vskip0truecm\noindent\font\BF=cmbx10 scaled \magstep2\BF #1\par}
\vskip1truecm
\leftline{\font\CSC=cmcsc10
{\CSC Friedrich Knop}}
\leftline{Department of Mathematics, Rutgers University, New Brunswick NJ
08903, USA}
\leftline{knop@math.rutgers.edu}
\vskip1truecm
\par}

\def\cite#1{\expandafter\ifx\csname#1\endcsname\relax
{\bf?}\immediate\write16{#1 ist nicht definiert!}\else\csname#1\endcsname\fi}
\def\expandwrite#1#2{\edef\next{\write#1{#2}}\next}
\def\neverexpand{\noexpand\noexpand\noexpand}
\def\strip#1\ {}
\def\ncite#1{\expandafter\ifx\csname#1\endcsname\relax
{\bf?}\immediate\write16{#1 ist nicht definiert!}\else
\expandafter\expandafter\expandafter\strip\csname#1\endcsname\fi}
\newwrite\AUX
\immediate\openout\AUX=\jobname.aux
\font\eightrm=cmr8\font\sixrm=cmr6
\font\eighti=cmmi8
\font\eightit=cmti8
\font\eightbf=cmbx8
\font\eightcsc=cmcsc10 scaled 833
\def\eightpoint{%
\textfont0=\eightrm\scriptfont0=\sixrm\def\rm{\fam0\eightrm}%
\textfont1=\eighti
\textfont\bffam=\eightbf\def\bf{\fam\bffam\eightbf}%
\textfont\itfam=\eightit\def\it{\fam\itfam\eightit}%
\def\csc{\eightcsc}%
\setbox\strutbox=\hbox{\vrule height7pt depth2pt width0pt}%
\normalbaselineskip=0,8\normalbaselineskip\normalbaselines\rm}
\newcount\absFnno\absFnno1
\write\AUX{\relFnno1}
\newif\ifMARKE\MARKEtrue
{\catcode`\@=11
\gdef\footnote{\ifMARKE\edef\@sf{\spacefactor\the\spacefactor}\/%
$^{\cite{Fn\the\absFnno}}$\@sf\fi
\MARKEtrue
\insert\footins\bgroup\eightpoint
\interlinepenalty100\let\par=\endgraf
\leftskip=0pt\rightskip=0pt
\splittopskip=10pt plus 1pt minus 1pt \floatingpenalty=20000\smallskip
\item{$^{\cite{Fn\the\absFnno}}$}%
\expandwrite\AUX{\neverexpand\ref{Fn\the\absFnno}{\neverexpand\the\relFnno}}%
\global\advance\absFnno1\write\AUX{\advance\relFnno1}%
\bgroup\strut\aftergroup\@foot\let\next}}
\skip\footins=12pt plus 2pt minus 4pt
\dimen\footins=30pc
\output={\plainoutput\immediate\write\AUX{\relFnno1}}
\newcount\Abschnitt\Abschnitt0
\def\beginsection#1. #2 \par{\advance\Abschnitt1%
\vskip0pt plus.10\vsize\penalty-250
\vskip0pt plus-.10\vsize\bigskip\vskip\parskip
\edef\TEST{\number\Abschnitt}
\expandafter\ifx\csname#1\endcsname\TEST\relax\else
\immediate\write16{#1 hat sich geaendert!}\fi
\expandwrite\AUX{\neverexpand\ref{#1}{\TEST}}
\leftline{\marginnote{#1}\bf\number\Abschnitt. \ignorespaces#2}%
\nobreak\smallskip\noindent\SATZ1\GNo0}
\def\Proof:{\par\noindent{\it Proof:}}
\def\Remark:{\ifdim\lastskip<\medskipamount\removelastskip\medskip\fi
\noindent{\bf Remark:}}
\def\Remarks:{\ifdim\lastskip<\medskipamount\removelastskip\medskip\fi
\noindent{\bf Remarks:}}
\def\Definition:{\ifdim\lastskip<\medskipamount\removelastskip\medskip\fi
\noindent{\bf Definition:}}
\def\Example:{\ifdim\lastskip<\medskipamount\removelastskip\medskip\fi
\noindent{\bf Example:}}
\def\Examples:{\ifdim\lastskip<\medskipamount\removelastskip\medskip\fi
\noindent{\bf Examples:}}
\newif\ifmarginalnotes\marginalnotesfalse
\newif\ifmarginalwarnings\marginalwarningstrue

\def\marginnote#1{\ifmarginalnotes\hbox to 0pt{\eightpoint\hss #1\ }\fi}

\def\strutdepth{\dp\strutbox}
\def\Randbem#1#2{\ifmarginalwarnings
{#1}\strut
\setbox0=\vtop{\eightpoint
\rightskip=0pt plus 6mm\hfuzz=3pt\hsize=16mm\noindent\leavevmode#2}%
\vadjust{\kern-\strutdepth
\vtop to \strutdepth{\kern-\ht0
\hbox to \hsize{\kern-16mm\kern-6pt\box0\kern6pt\hfill}\vss}}\fi}

\def\Zitat!{\Randbem{\bf?}{\bf Zitat}}

\newcount\SATZ\SATZ1
\def\proclaim #1. #2\par{\ifdim\lastskip<\medskipamount\removelastskip
\medskip\fi
\noindent{\bf#1.\ }{\it#2\Par}
\ifdim\lastskip<\medskipamount\removelastskip\goodbreak\medskip\fi}
\def\Aussage#1{\expandafter\def\csname#1\endcsname##1.{\resetitem
\ifx?##1?\relax\else
\edef\TEST{#1\penalty10000\ \number\Abschnitt.\number\SATZ}
\expandafter\ifx\csname##1\endcsname\TEST\relax\else
\immediate\write16{##1 hat sich geaendert!}\fi
\expandwrite\AUX{\neverexpand\ref{##1}{\TEST}}\fi
\proclaim {\marginnote{##1}\number\Abschnitt.\number\SATZ. #1\global\advance\SATZ1}.}}
\Aussage{Theorem}
\Aussage{Proposition}
\Aussage{Corollary}
\Aussage{Lemma}
\font\la=lasy10
\def\strich{\hbox{$\vcenter{\hbox
to 1pt{\leaders\hrule height -0,2pt depth 0,6pt\hfil}}$}}
\def\dashedrightarrow{\hbox{%
\hbox to 0,5cm{\leaders\hbox to 2pt{\hfil\strich\hfil}\hfil}%
\kern-2pt\hbox{\la\char\string"29}}}

\def\Bindestrich{\penalty10000-\hskip0pt}
\let\_=\Bindestrich
\def\.{{\sfcode`.=1000.}}

\def\Par{\par}
\def\:={\mathrel{\raise0,9pt\hbox{.}\kern-2,77779pt
\raise3pt\hbox{.}\kern-2,5pt=}}
\def\=:{\mathrel{=\kern-2,5pt\raise0,9pt\hbox{.}\kern-2,77779pt
\raise3pt\hbox{.}}} 
\def\into{\hookrightarrow}
\def\pfeil{\rightarrow}

\def\Ugleich{\hbox{$\cup$\kern.5pt\vrule depth -0.5pt}}
\def\|#1|{\mathop{\rm#1}\nolimits}
\def\<{\langle}
\def\>{\rangle}
\let\Times=\times
\def\times{\mathop{\Times}}
\let\Otimes=\otimes
\def\otimes{\mathop{\Otimes}}
\catcode`\@=11
\def\hex#1{\ifcase#1 0\or1\or2\or3\or4\or5\or6\or7\or8\or9\or A\or B\or
C\or D\or E\or F\else\message{Warnung: Setze hex#1=0}0\fi}
\def\fontdef#1:#2,#3,#4.{%
\alloc@8\fam\chardef\sixt@@n\FAM
\ifx!#2!\else\expandafter\font\csname text#1\endcsname=#2
\textfont\the\FAM=\csname text#1\endcsname\fi
\ifx!#3!\else\expandafter\font\csname script#1\endcsname=#3
\scriptfont\the\FAM=\csname script#1\endcsname\fi
\ifx!#4!\else\expandafter\font\csname scriptscript#1\endcsname=#4
\scriptscriptfont\the\FAM=\csname scriptscript#1\endcsname\fi
\expandafter\edef\csname #1\endcsname{\fam\the\FAM\csname text#1\endcsname}
\expandafter\edef\csname hex#1fam\endcsname{\hex\FAM}}
\catcode`\@=12 

\fontdef Ss:cmss10,,.
\fontdef Fr:eufm10,eufm7,eufm5.


%
\fontdef bbb:msbm10,msbm7,msbm5.
\fontdef mbf:cmmib10,cmmib7,.
\fontdef msa:msam10,msam7,msam5.

\def\FF{{\bbb F}}

\def\NN{{\bbb N}}

\def\ZZ{{\bbb Z}}
\def\cA{{\cal A}}

\def\cN{{\cal N}}
\def\cS{{\cal S}}\def\cT{{\cal T}}

\mathchardef\leer=\string"0\hexbbbfam3F
\mathchardef\subsetneq=\string"3\hexbbbfam24
\mathchardef\semidir=\string"2\hexbbbfam6E
\mathchardef\dirsemi=\string"2\hexbbbfam6F
\mathchardef\haken=\string"2\hexmsafam78
\mathchardef\auf=\string"3\hexmsafam10
\let\OL=\overline
\def\overline#1{{\hskip1pt\OL{\hskip-1pt#1\hskip-.3pt}\hskip.3pt}}


%
\newdimen\Parindent
\Parindent=\parindent


\abovedisplayskip 9.0pt plus 3.0pt minus 3.0pt
\belowdisplayskip 9.0pt plus 3.0pt minus 3.0pt
\newdimen\Grenze\Grenze2\Parindent\advance\Grenze1em
\newdimen\Breite
\newbox\DpBox
\def\NewDisplay#1
#2$${\Breite\hsize\advance\Breite-\hangindent
\setbox\DpBox=\hbox{\hskip2\Parindent$\displaystyle{%
\ifx0#1\relax\else\eqno{#1}\fi#2}$}%
\ifnum\predisplaysize<\Grenze\abovedisplayskip\abovedisplayshortskip
\belowdisplayskip\belowdisplayshortskip\fi
\global\futurelet\nexttok\WEITER}
\def\WEITER{\ifx\nexttok\qed\expandafter\leftQEDdisplay
\else\leftdisplay\fi}
\def\leftdisplay{\hskip-\hangindent\leftline{\box\DpBox}$$}
\def\leftQEDdisplay{\hskip-\hangindent
\line{\copy\DpBox\hfill\lower\dp\DpBox\copy\QEDbox}%
\belowdisplayskip0pt$$\bigskip\let\nexttok=}
\everydisplay{\NewDisplay}
\newcount\GNo\GNo=0
\newcount\maxEqNo\maxEqNo=0
\def\eqno#1{%
\global\advance\GNo1
\edef\FTEST{$(\number\Abschnitt.\number\GNo)$}
\ifx?#1?\relax\else
\ifnum#1>\maxEqNo\global\maxEqNo=#1\fi%
\expandafter\ifx\csname E#1\endcsname\FTEST\relax\else
\immediate\write16{E#1 hat sich geaendert!}\fi
\expandwrite\AUX{\neverexpand\ref{E#1}{\FTEST}}\fi
\llap{\hbox to 40pt{\ifx?#1?\relax\else\marginnote{E#1}\fi\FTEST\hfill}}}

\catcode`@=11
\def\eqalignno#1{\null\!\!\vcenter{\openup\jot\m@th\ialign{\eqno{##}\hfil
&\strut\hfil$\displaystyle{##}$&$\displaystyle{{}##}$\hfil\crcr#1\crcr}}\,}
\catcode`@=12

\newbox\QEDbox
\newbox\nichts\setbox\nichts=\vbox{}\wd\nichts=2mm\ht\nichts=2mm
\setbox\QEDbox=\hbox{\vrule\vbox{\hrule\copy\nichts\hrule}\vrule}
\def\qed{\leavevmode\unskip\hfil\null\nobreak\hfill\copy\QEDbox\medbreak}
\newdimen\HIindent
\newbox\HIbox
\def\setHI#1{\setbox\HIbox=\hbox{#1}\HIindent=\wd\HIbox}
\def\HI#1{\par\hangindent\HIindent\hangafter=0\noindent\leavevmode
\llap{\hbox to\HIindent{#1\hfil}}\ignorespaces}

\newdimen\maxSpalbr
\newdimen\altSpalbr
\newcount\Zaehler


\newif\ifxxx

{\catcode`/=\active

\gdef\beginrefs{%
\xxxfalse
\catcode`/=\active
\def/{\string/\ifxxx\hskip0pt\fi}
\def\TText##1{{\xxxtrue\tt##1}}
\expandafter\ifx\csname Spaltenbreite\endcsname\relax
\def\Spaltenbreite{1cm}\immediate\write16{Spaltenbreite undefiniert!}\fi
\expandafter\altSpalbr\Spaltenbreite
\maxSpalbr0pt
\gdef\alt{}
\def\\##1\relax{%
\gdef\neu{##1}\ifx\alt\neu\global\advance\Zaehler1\else
\xdef\alt{\neu}\global\Zaehler=1\fi\xdef\SigText{##1\the\Zaehler}}
\def\L|Abk:##1|Sig:##2|Au:##3|Tit:##4|Zs:##5|Bd:##6|S:##7|J:##8|xxx:##9||{%
\def\SigText{##2}\global\setbox0=\hbox{##2\relax}
\edef\TEST{[\SigText]}
\expandafter\ifx\csname##1\endcsname\TEST\relax\else
\immediate\write16{##1 hat sich geaendert!}\fi
\expandwrite\AUX{\neverexpand\ref{##1}{\TEST}}
\setHI{[\SigText]\ }
\ifnum\HIindent>\maxSpalbr\maxSpalbr\HIindent\fi
\ifnum\HIindent<\altSpalbr\HIindent\altSpalbr\fi
\HI{\marginnote{##1}[\SigText]}
\ifx-##3\relax\else{##3}: \fi
\ifx-##4\relax\else{##4}{\sfcode`.=3000.} \fi
\ifx-##5\relax\else{\it ##5\/} \fi
\ifx-##6\relax\else{\bf ##6} \fi
\ifx-##8\relax\else({##8})\fi
\ifx-##7\relax\else, {##7}\fi
\ifx-##9\relax\else, \TText{##9}\fi\Par}
\def\B|Abk:##1|Sig:##2|Au:##3|Tit:##4|Reihe:##5|Verlag:##6|Ort:##7|J:##8|xxx:##9||{%
\def\SigText{##2}\global\setbox0=\hbox{##2\relax}
\edef\TEST{[\SigText]}
\expandafter\ifx\csname##1\endcsname\TEST\relax\else
\immediate\write16{##1 hat sich geaendert!}\fi
\expandwrite\AUX{\neverexpand\ref{##1}{\TEST}}
\setHI{[\SigText]\ }
\ifnum\HIindent>\maxSpalbr\maxSpalbr\HIindent\fi
\ifnum\HIindent<\altSpalbr\HIindent\altSpalbr\fi
\HI{\marginnote{##1}[\SigText]}
\ifx-##3\relax\else{##3}: \fi
\ifx-##4\relax\else{##4}{\sfcode`.=3000.} \fi
\ifx-##5\relax\else{(##5)} \fi
\ifx-##7\relax\else{##7:} \fi
\ifx-##6\relax\else{##6}\fi
\ifx-##8\relax\else{ ##8}\fi
\ifx-##9\relax\else, \TText{##9}\fi\Par}
\def\Pr|Abk:##1|Sig:##2|Au:##3|Artikel:##4|Titel:##5|Hgr:##6|Reihe:{%
\def\SigText{##2}\global\setbox0=\hbox{##2\relax}
\edef\TEST{[\SigText]}
\expandafter\ifx\csname##1\endcsname\TEST\relax\else
\immediate\write16{##1 hat sich geaendert!}\fi
\expandwrite\AUX{\neverexpand\ref{##1}{\TEST}}
\setHI{[\SigText]\ }
\ifnum\HIindent>\maxSpalbr\maxSpalbr\HIindent\fi
\ifnum\HIindent<\altSpalbr\HIindent\altSpalbr\fi
\HI{\marginnote{##1}[\SigText]}
\ifx-##3\relax\else{##3}: \fi
\ifx-##4\relax\else{##4}{\sfcode`.=3000.} \fi
\ifx-##5\relax\else{In: \it ##5}. \fi
\ifx-##6\relax\else{(##6)} \fi\PrII}
\def\PrII##1|Bd:##2|Verlag:##3|Ort:##4|S:##5|J:##6|xxx:##7||{%
\ifx-##1\relax\else{##1} \fi
\ifx-##2\relax\else{\bf ##2}, \fi
\ifx-##4\relax\else{##4:} \fi
\ifx-##3\relax\else{##3} \fi
\ifx-##6\relax\else{##6}\fi
\ifx-##5\relax\else{, ##5}\fi
\ifx-##7\relax\else, \TText{##7}\fi\Par}
\bgroup
\baselineskip12pt
\parskip2.5pt plus 1pt
\hyphenation{Hei-del-berg Sprin-ger}
\sfcode`.=1000
\beginsection References. References

}}

\def\endrefs{%
\expandwrite\AUX{\neverexpand\ref{Spaltenbreite}{\the\maxSpalbr}}
\ifnum\maxSpalbr=\altSpalbr\relax\else
\immediate\write16{Spaltenbreite hat sich geaendert!}\fi
\egroup\write16{Letzte Gleichung: E\the\maxEqNo}
\write16{Letzte Aufzaehlung: I\the\maxItemcount}}



\def\1{{\bf 1}}
\def\Ah{{\hat\cA}}
\def\cT{{\cal T}}
\def\Ap{{\cA_p}}

\title{A construction of semisimple tensor categories}
\bigskip

{\narrower
\noindent{\bf Abstract.}
Let $\cA$ be an abelian category such that every object has
only finitely many subobjects. From $\cA$ we construct a semisimple tensor
category $\cT$. We show that $\cT$ interpolates the categories
$\|Rep|(\|Aut|(p),K)$ where $p$ runs through certain projective
pro\_objects of $\cA$. This extends a construction of Deligne
for symmetric groups.

}

\beginsection intro. Introduction

Let $K$ be a field of characteristic zero. In \cite{De}, Deligne
constructed a tensor category $\|Rep|(S_t,K)$ over $K$ depending on a
parameter $t\in K$. If $t\not\in\NN$ then $\|Rep|(S_t,K)$ is
abelian semisimple. Otherwise, it has as quotient the category of finite
dimensional representations of the symmetric group $S_t$.

In this paper, we extend Deligne's construction. Starting from an
abelian category $\cA$ such that every object has only finitely many
subobjects we construct a tensor category $\cT=\cT(\cA,K)$ which
depends on parameters $t_\phi$, one for each isomorphism class of
simple objects in $\cA$. We show that $\cT$ is semisimple if none of
the parameters is singular (see \S3 for a definition). Then the
simple objects of $\cT$ correspond to pairs $(x,\pi)$ where $x$ is an
object of $\cA$ and $\pi$ is an irreducible representation of
$\|Aut|_\cA(x)$. If all parameters are singular, the category $\cT$
has as quotient $\|Rep|(\|Aut|(p),K)$ where $p$ is a projective
(pro-)object of $\cA$.

The main example is $\cA=\|Mod|(\FF_q)$, the category of finite
dimensional $\FF_q$\_vector spaces. In that case, the simple objects
of $\cT$ correspond to irreducible representations of $GL(m,\FF_q)$,
$m\in\NN$. There is only one parameter and this parameter is singular
if and only it is a power $q^n$. In that case, the category
has as a quotient $\|Rep|(GL(n,\FF_q),K)$. This proves conjecture
\cite{De}, p. 3, (A) of Deligne. In unpublished work, Deligne has
also proved his conjecture.

\beginsection CoCor. The construction of $\cT(\cA,K)$

Let $\cA$ be an essentially small abelian category such that every
object has finite length. Let $\Ah$ be the set of isomorphism classes
of simple objects of $\cA$. The category we are going to construct
will depend on $\cA$ and on a fixed map $\Ah\pfeil K$ where $K$ is a
field. The image of $\phi\in\Ah$ in $K$ will be denoted by $t_\phi$. Let
$\kappa(\cA)$ be the free commutative monoid generated by $\Ah$. It
coincides with the Grothendieck monoid of $\cA$. Thus, every object
$x$ gives rise to an element $\<x\>$ of $\kappa(\cA)$. In particular,
$\<x\>=\<y\>$ if and only if $x$ and $y$ have the same composition
factors.

A {\it correspondence} between two objects $x$ and $y$ is a morphism
$F:c\pfeil x\oplus y$. If $F$ is a monomorphism then it is called a
{\it relation}. Two correspondences $F:c\pfeil x\oplus y$ and
$G:d\pfeil x\oplus y$ are called {\it equivalent} if
$$
\|im|F=\|im|G\quad{\rm and}\quad\<\|ker|F\>=\<\|ker|G\>
$$
We define the pseudo-abelian tensor category $\cT$ in several steps. First
we define the category $\cT_0$:

\smallskip\hangindent20pt\hangafter1\noindent{\it Objects:} Same as
$\cA$. The object $x$ of $\cA$, regarded as an object of $\cT_0$, will
be denoted by $[x]$.

\hangindent20pt\hangafter1\noindent{\it Morphisms:}
Equivalence classes of correspondences.

\hangindent20pt\hangafter1\noindent{\it Composition:} If
$G:c\pfeil x\oplus y$ and $F:c\pfeil y\oplus z$ are correspondences
then $FG$ is the equivalence class of $c\times_y d\pfeil
x\oplus z$. It is easy to see that the composition is well defined and
associative.

\smallskip

The category $\cT_0$ becomes a symmetric monoidal category by defining
$[x]\otimes[y]:=[x\oplus y]$. The unit object is $\1=[0]$.  Each
object is selfdual with $\delta:\1\pfeil[x]\otimes[x]$ and
$\|ev|:[x]\otimes[x]\pfeil\1$ given by the diagonal morphism $x\pfeil
x\oplus x$.

Let $\cT_1(\cA)$ be the category with the same objects as $\cT_0$ but
with $\|Hom|_{\cT_1(\cA)}([x],[y])$ being the free abelian group
generated by $\|Hom|_{\cT_0}([x],[y])$. The ring
$\|End|_{\cT_1(\cA)}(\1)$ is isomorphic to the polynomial ring
$\ZZ[\Ah]$. Thus, every $\|Hom|$-space is a $\ZZ[\Ah]$\_module. The
fixed map $\Ah\pfeil K$ induces a homomorphism $\ZZ[\Ah]\pfeil K$. Now
we define the category $\cT_1(\cA,K)$ as having the same objects as
$\cT_1(\cA)$ but with morphisms
$$
\|Hom|_{\cT_1(\cA,K)}([x],[y])=
\|Hom|_{\cT_1(\cA)}([x],[y])\otimes_{\ZZ[\Ah]}K.
$$

Finally, let $\cT=\cT(\cA,K)$ be the pseudo\_abelian completion of
$\cT_1(\cA,K)$, i.e., the category obtained by adjoining finite direct
sums and images of idempotents. The tensor product on $\cT_0$ induces
a symmetric $K$\_bilinear tensor product on $\cT$ such that every
object has a dual.

\beginsection MainSec. The semisimplicity of $\cT(\cA,K)$ for regular
parameters

We call $\cA$ {\it finitary} if every object has only finitely many
subobjects. In that case, all $\|Hom|$\_spaces of $\cA$ are finite and
all $\|Hom|$\_spaces of $\cT(\cA,K)$ are
finite dimensional $K$\_vector spaces. For $\phi\in\Ah$ let
$m_\phi$ be a simple object in the isomorphism class $\phi$. Then
$\|End|_\cA(m_\phi)$ is a finite field. Let $q_\phi$ be its order.

\Definition:
An element $t\in K$ is {\it $\phi$\_singular} if
either
\Item{} $t\in\{1,q_\phi,q_\phi^2,\ldots\}$ or
\Item{} $t=0$ and $m_\phi$ is part of a non\_splitting short exact
sequence.

\Theorem main.
Let $\cA$ be an essentially small finitary abelian category. Assume
that $K$ is a field of characteristic zero and that there is no
$\phi\in\Ah$ such that $t_\phi$ is $\phi$\_singular.
Then $\cT(\cA,K)$ is a semisimple
tensor category. The simple objects correspond to isomorphism classes
of pairs $(x,\pi)$ where $x$ is an object of $\cA$ and $\pi$ is an
irreducible representation (over $K$) of $\|Aut|_\cA(x)$.

The proof has two main ingredients.

\Lemma ingred1.
The pairing $\|Hom|_\cT(\1,X)\times\|Hom|_\cT(X,\1)\pfeil K:
(F,G)\mapsto\|tr|GF$ is perfect for all $X\in\|Ob|\cT$.

\Proof: It suffices to check this for $X=[x]$. Then the assertion
boils down to the non\_vanishing of the determinant
$\Delta_x:=\|det|\big(\<u\cap v\>_K\big)_{u,v\subseteq x}$.  Here
$\<u\>_K$ denotes the image of $\<u\>\in\kappa(\cA)\subseteq\ZZ[\Ah]$
in $K$. A formula of Lindstr\"om \cite{Li} and Wilf \cite{Wi} implies
$$50
\Delta_x=\prod\nolimits_{y\subseteq x}p_y\quad\hbox{with}\quad
p_y:=\sum\nolimits_{u\subseteq y}\mu(u,y)\<u\>_K.
$$
Here $\mu(u,y)$ is the M\"obius function of the subobject lattice of
$x$ (or $y$). Let $m$ be a simple subobject of $y$. Then $p_y$
factorizes as $p_y=(t_\phi-\alpha)\,p_{y/m}$ where $\alpha$ is the
number of complements of $m$ in $y$ (Stanley \cite{St1}, \cite{St2}). Since
$\alpha$ is either $0$ or a power of $q_\phi$ we conclude by
induction.\qed

Now we come to the second main ingredient. Let $\ell(x)$ denote the length
of an object $x$ of $\cA$.

\Definition:
For a $\cT$\_morphism $F$ let $\ell(F)$ be the least number $l$
such that $F$ factorizes through $[x_1]\oplus\ldots\oplus[x_s]$ with
$\ell(x_i)\le l$ for all $i$.

\medskip

For the next statement, note that $x\mapsto[x]$,
$f\mapsto\|graph|(f)$ defines an embedding of $\cA$ into $\cT$.

\Lemma ingred2.
Let $x$ and $y$ be two objects of $\cA$ with
$\ell(x)=\ell(y)=l$. Then
$$40
\|Hom|_\cT([x],[y])=K[\|Isom|_\cA(x,y)]
\oplus\{F\in\|Hom|_\cT([x],[y])\mid\ell(F)<l\}.
$$

\Proof: For a correspondence $F:c\pfeil x\oplus y$ with components
$F_x:c\pfeil x$ and $F_y:c\pfeil y$ let
$$
\|core|F:=c/(\|ker|F_x+\|ker|F_y).
$$
Now \cite{E40} follows from the following claims:

\smallskip

\Item{6}$F$ factorizes in $\cT_0$ through $[\|core|F]$.

\Item{7}If $F$ factorizes in $\cT_0$ through $[z]$ then
$\ell(z)\ge\ell(\|core|F)$.

\Item{8}$F\in\|Isom|_\cA(x,y)$ if and only if $\ell(\|core|F)=l$.\qed

With these two propositions at hand, the proof of \cite{main}
proceeds along the same lines as that of \cite{De} Thm.~2.18.\qed

\beginsection specialization. Specialization of $\cT(\cA,K)$ at singular parameters

In this section we study the category $\cT$ when all parameters
$t_\phi$ are singular but not zero. More precisely assume
$t_\phi=q_\phi^{r_\phi}\hbox{ with }r_\phi\in\NN\hbox{ for all
}\phi\in\Ah$.  Let $\Ap$ be the category of all pro-objects of
$\cA$. This category has enough projectives. Let
$$28
p\auf\bigoplus\nolimits_{\phi\in\Ah}m_\phi^{\oplus r_\phi}
$$
be a projective cover and put $A(p):=\|Aut|_\Ap(p)$, a profinite
group. Let $\cN(Y,X)$ be the set of all
$\cT$\_morphisms $F:Y\pfeil X$ with $\|tr|GF=0$ for all $G:X\pfeil Y$.
Then $\cN$ forms a tensor ideal of $\cT$.

\Theorem specializ.
There is a functor $\cS:\cT\pfeil\|Rep|(A(p),K)$ which
identifies $\cT/\cN$ with $\|Rep|(A(p),K)$.

\Proof: For a set $S$ let $K[S]$ be the space of all functions $S\pfeil
K$. First we define a functor $\cS:\cT_0\pfeil\|Rep|(A(p),K)$. Let $x$ be
an object of $\cA$ and $F:c\pfeil x\oplus y$ a relation. Then we define
$$
\cS([x]):=K[\|Hom|_\Ap(p,x)],\quad
\cS(F):\cS([x])\pfeil\cS([y]):\alpha\mapsto \sum_{{\beta:p\pfeil
c\atop\alpha=F_x\beta}}(F_y\beta)
$$
Using the projectivity of $p$ one checks that $\cS$ is a well defined
tensor functor. Next, \cite{E28} implies that
$\cS(x\pfeil0\oplus0)=\<x\>_K\in K$. This ensures that $\cS$ extends
uniquely to $\cT_1$ and then to $\cT$. That $\cS$ has the stated
property follows along the same lines as the proof of
\cite{De} Thm.~6.2.\qed

\beginsection ExRem. Examples and final remarks

As for examples, we already mentioned $\cA=\|Mod|(\FF_q)$
in the introduction. We point out three more:

\smallskip\noindent 1. Let $\cA$ be the category of homomorphisms $U\pfeil V$
between $\FF_q$\_vector spaces. Then $\cT$ interpolates the
representations of the parabolic $\pmatrix{*&*\cr0&*\cr}\subseteq
GL(n_1+n_2,\FF_q)$ with arbitrary block sizes $n_1$ and $n_2$. The set
$\Ah$ consists of two elements. The regular parameters are
$t_1,t_2\ne0,1,q,q^2,\ldots$.

\smallskip\noindent 2. Let $\cA$ be the category of pairs $(V,\alpha)$ where $V$
is an $\FF_q$\_vector space and $\alpha$ is a nilpotent endomorphism
of $V$. Then $\cT$ interpolates $\|Rep|(GL(n,\FF_q[\![x]\!]),K)$ with
regular values $t\ne0,1,q,q^2,\ldots$.

\smallskip\noindent 3. Let $\cA$ be the category of all finite abelian
$p$\_groups. Then $\cT$ interpolates $\|Rep|(GL(n,\hat\ZZ_p),K)$ with
regular values $t\ne0,1,p,p^2,\ldots$.

\Remarks:
1. Deligne's category $\|Rep|(S_t,K)$ is obtained by taking
for $\cA$ the opposite of the category of finite sets. Of course, this
category is not abelian. However, most of the results work more
generally in the framework of exact Mal'cev categories in the sense of
\cite{CLP}. These comprise not only all abelian categories and the
opposite of the category of sets but also the categories of finite
groups, finite rings and many more. In particular, it is possible to
interpolate the representation categories of wreath products $S_n\wr
G$ (for fixed $G$) or $S_{n_1}\wr S_{n_2}\wr S_{n_3}\ldots$ Details
will appear elsewhere.

\smallskip\noindent
2. The basic objects in \cite{De} are slightly different. Let $x$ be an
object of $\cA$. Then every subobject $y$ of $x$ gives rise to an
idempotent via the relation $y\into x\oplus x$. These idempotents
commute and induce a decomposition $[x]=\mathop\oplus_{y\subseteq
x}[y]^*$. It is these $[y]^*$, Deligne is working with. Observe that
$\cS([x]^*)$ is the space of functions on the set of all {\it epimorphisms}
$p\auf x$. The factorization \cite{E50} is an easy
consequence of the decomposition.

\beginrefs

\L|Abk:CLP|Sig:CLP|Au:Carboni, A.; Lambek, J.; Pedicchio,
M.|Tit:Diagram chasing in Mal'cev categories|Zs:J. Pure Appl.
Algebra|Bd:69|S:271--284|J:1991|xxx:-||

\L|Abk:De|Sig:De|Au:Deligne, P.|Tit:La cat\'egorie des
repr\'esentations du groupe sym\'etrique $S_t$ lorsque $t$ n'est pas
un entier naturel|Zs:Preprint|Bd:-|S:78
pages|J:-|xxx:www.math.ias.edu/\lower4pt\hbox{\char126}phares/deligne/Symetrique.pdf||

\L|Abk:Li|Sig:Li|Au:Lindstr\"om, B.|Tit:Determinants on
semilattices|Zs:Proc. Amer Math. Soc.|Bd:20|S:207--208|J:1969|xxx:-||

\L|Abk:St1|Sig:\\St|Au:Stanley, R.|Tit:Modular elements of geometric
lattices|Zs:Algebra Universalis|Bd:1|S:214-217|J:1971/72|xxx:-||

\L|Abk:St2|Sig:\\St|Au:Stanley, R.|Tit:Supersolvable
lattices|Zs:Algebra Universalis|Bd:2|S:197--217|J:1972|xxx:-||

\L|Abk:Wi|Sig:Wi|Au:Wilf, H.|Tit:Hadamard determinants, M\"obius
functions, and the chromatic number of a
graph|Zs:Bull. Amer. Math. Soc.|Bd:74|S:960--964|J:1968|xxx:-||

\endrefs

\bye